\documentstyle[12pt]{article}
\def\IR{\hbox{\rm I\kern-.18em R}}

\def\IC{\hbox{\rm C\kern-.43em
       \vrule depth 0ex height 1.4ex width .05em\kern.41em}}

\def\IF{\hbox{\rm I\kern-.18em F}} %symbol for field

\def\IN{\hbox{\rm I\kern-.17em N}}

\def\ZZ{\hbox{\rm Z\kern-.3em Z}}
\def\IQ{\hbox{\rm Q\kern-.43em
       \vrule depth 0ex height 1.4ex width .05em\kern.41em}}
\def\hh{\hrule height0.9pt width1.1em}
\def\vv{\vrule width0.8pt depth0.12em height0.92em}
\def\square{\vbox{\kern0.15em\hh\kern0.9em\hh\kern-1.05em
            \hbox{\vv\kern0.93em\vv}}}
\def\qed{{\hfill\square}}
\def\dim{{\rm dim}\,}
\def\span{{\rm span}\,}
\def\supp{{\rm supp}\,}
\def\card{{\rm card}\,}
\newcommand{\lra}{\longrightarrow}
\newcommand{\w}{\widetilde }

\newcommand{\bbF}{\IF}

\newcommand{\bbN}{\IN}
\newcommand{\bbR}{\IR}

\begin{document}
%%%%%%%%%%%%%%%%%%%%%%%%%%%%%%%%%%%%
%%%%	BEGIN DOCUMENT HERE	%%%%
%%%%%%%%%%%%%%%%%%%%%%%%%%%%%%%%%%%%
\parskip=0.1in
\baselineskip 15pt
\centerline{\bf ISOMETRIES OF CROSS PRODUCTS OF SEQUENCE SPACES}

\medskip
\centerline{\sc Chi-Kwong Li \quad and \quad Beata Randrianantoanina}

\medskip
%\centerline{\today}

\medskip\noindent
{\bf Abstract}

\narrower \smallskip  \sl
Let $X_0, X_1, \dots, X_k$ with $k \in \IN\cup\{\infty\}$ be sequence spaces 
$($finite or infinite dimensional$)$ over $\IC$ or $\IR$ with absolute norms 
$N_i$ for $i = 0, \dots, k$, $($i.e., with  1-unconditional bases$)$ such 
that $\dim X_0 = k$. Define an absolute norm on the  cross product space 
$($also known as 
the $X_0$ 1-unconditional sum$)$ $X_1 \times \cdots \times X_k$ by 
$$
N(x_1, \dots, x_k) = N_0(N_1(x_1), \dots, N_k(x_k)) \quad \hbox{ for  
all } \quad (x_1, \dots, x_k) \in X_1 \times \cdots \times X_k.
$$ 
We show that every sequence space with an absolute norm has an intrinsic
cross product structure of this form.  The result is used to prove a
characterization of isometries of complex cross product spaces that
covers all the existing results.  We demonstrate by examples and the 
theory of finite reflection groups that it is impossible to extend the 
complex result to the real case.  Nevertheless, some new isometry theorems
are obtained for  real cross product spaces. \smallskip
\rm 



\leftskip=0pt\rightskip=0pt
\bigskip\noindent
{\bf 1. Introduction}
\rm

Let $X_0, X_1, \dots, X_k$ with $k \in \IN\cup\{\infty\}$ be sequence  
spaces (finite or infinite dimensional) over $\IC$ or $\IR$ with absolute 
norms $N_i$ for $i = 0, \dots, k$, (i.e., with 1-unconditional bases)  
such that $\dim X_0 = k$. Define an absolute norm on the cross product 
space (also known as the $X_0$ 1-unconditional sum)  
$X_1 \times \cdots \times X_k$ by 
$$
N(x_1, \dots, x_k) = N_0(N_1(x_1), \dots, N_k(x_k)) \quad \hbox{ for  
all} 
\quad (x_1, \dots, x_k) \in X_1 \times \cdots \times X_k.
$$ 
Such a space is denoted by $X_0(X_1, \dots, X_k)$.
If $X_1 = \cdots = X_k = Y$, the notation $X_0(Y)$ is used.
The purpose of this paper is to study the geometry and isometries 
of $X_0(X_1, \dots, X_k)$.

The study of cross products of normed spaces arises naturally in many  
areas of mathematics. In particular, they have been a source of examples and 
counter-examples in geometric theory of Banach spaces (see e.g. 
\cite{Day41,DV86}).  

To understand the geometry of a normed vector space, it is useful to know the 
structure of its isometries.  In fact, many authors have studied the 
isometries of cross product spaces. For example, Fleming, Goldstein,  Jamison 
\cite{FGJ} studied isometries of 1-unconditional sums of Hilbert  spaces (see 
also Fleming and Jamison \cite{FJ74a,FJ74}) in the complex case and  Rosenthal 
\cite{R86} obtained the result for the real case, Greim \cite{Gr}  studied 
surjective isometries of $\ell_p$ sums of Banach spaces (see also  
\cite{KL}), Fleming and Jamison \cite{FJ89} studied isometries of complex  
$c_0-$sums and $E-$sums, where $E$ is ``sufficiently $\ell_p$ like'',  say, 
$E$ is a  ``nice'' Orlicz space (see \cite{FJ89} for precise definitions).  It 
turns out that all the results in these papers show that a surjective isometry 
always preserves the cross product structure of the space. There is also a 
number of papers that address this problem in non-atomic function spaces. For 
the detailed  discussion of the literature we refer the readers to the survey 
\cite{FJs}. 

Since $N_0, \dots, N_k$ are absolute norms, the norm $N$ on $X_0(X_1, \dots, 
X_k)$ is also absolute. In the very interesting paper of Schneider and Turner 
\cite{ST}, the  authors determine the structure of  isometries for an absolute 
norm $N$ on  $\IC^n$, which is the space of complex column vectors with $n$ 
entries and  will be viewed as an $n$-dimensional sequence space in our 
discussion. In  particular, it was shown (cf. \cite[(2.3) and (7.7)]{ST}) that 
if the absolute norm is {\it normalized} so that $N(e_i) = 1$ for all standard 
unit vectors for $1 \le i \le  n$, then $\IC^n$ can be decomposed into a 
direct sum of $Y_i = \span\{v: v \in  E_i\}$ for $i = 1, \dots, k$, where 
$E_1 \cup \cdots \cup E_k = \{e_1,  \dots, e_n\}$, the standard basis of 
$\IC^n$, and there exists an absolute norm $N_0$ on $\IC^k$ such that 
\begin{itemize}
\item[(a)] 
each $(Y_i, N)$ is just an $\ell_2$ space, i.e., the Euclidean space, 
and
\item[(b)] $N(x_1,\dots, x_k) = N_0(N(x_1), \dots, N(x_k))$ for each
$x = (x_1, \dots, x_k) \in Y_1\times \cdots \times Y_k \cong \IC^n$.
\end{itemize}
Furthermore, an isometry for $N$ must be of the form
\begin{equation}\label{1}
(x_1, \dots, x_k) \mapsto (U_1 x_{\pi(1)}, \dots, U_k x_{\pi(k)})
\end{equation}
for some unitary $U_i$, $1 \le i \le k$, and a permutation $\pi$ of
the set $\{1, \dots, k\}$ such that
$N_0(z_1,\ldots,z_k) = N_0(z_{\pi (1)}, \ldots, z_{\pi (k)})$. 

\smallskip
This result was later extended to infinite dimensional spaces by 
Kalton and Wood \cite[Theorem 6.1]{KW}. 

\smallskip
By the above result, one sees that there is an intrinsic cross  product 
structure on every complex sequence  space with an absolute norm, 
and such a structure 
is useful in characterizing  isometries.  However, the  cross product 
decomposition in \cite{ST,KW} can only identify $\ell_2$ components.  If such 
components do not exist, then every {\it factor} (or  {\it factor space}) 
$X_i$ will be one dimensional, and the decomposition will not be very 
interesting. Of course, one can still get the very useful conclusion that 
every isometry for the norm must be a signed  permutation 
operator, i.e., an operator of the form (\ref{1}) with all $Y_i$  
being 1-dimensional vector spaces (scalars). Nevertheless, the theorem in 
\cite{ST} and \cite{KW} seems inadequate to explain the various isometry 
results on  cross product spaces. 

In this paper, we propose a new way to decompose a real or complex  sequence 
space with an absolute norm into a cross product of simpler spaces, which are 
not necessarily Euclidean.  Using this decomposition, we obtain a 
characterization  of the isometries  similar to that in \cite{ST,KW} in the 
complex case that covers all the known  isometry results on cross product 
spaces.  The decomposition also allows us to  obtain new and reprove existing 
isometry results on real cross product  spaces. 

Our paper is organized as follows. In Section 2,  we show that every  (complex 
or real) sequence space with an absolute norm can be decomposed into a cross 
product of simpler spaces, which could possibly be further presented as cross 
products of subsequent simpler spaces. Thus we get a ``multi--level''  cross
product structure present in every space (see Remark 1 after Theorem 2.4). 

In Section 3, we prove that in a complex sequence space with an  absolute 
norm, every surjective isometry necessarily preserves the intrinsic cross  
product structure described in  Theorem~2.4. A number of corollaries covering 
various existing isometry results on complex cross product spaces are also 
presented. 

In section 4,  we study  isometries of real spaces with cross product 
structure.   In particular, we show that our complex result can be applied to 
all real spaces whose isometry group is contained  in the group of signed 
permutations. This includes in particular spaces with $\Delta$-bases 
\cite{GL79} and spaces which are $p$-convex with  constant 1 for $2<p<\infty$ 
\cite{Rseq}. 

However the situation in real spaces is more complicated since there  are many 
natural spaces with cross product structure which have isometries  other than 
signed permutation operators i.e.,   isometries do not always preserve 
disjointness of vectors (see the examples in Section 4). 

Also, note that there are real spaces with explicit cross product structure 
which is not preserved by some isometries (see Examples 4  and 5 in Section 
4). We feel that such pathology should be rare, but  since  every finite group 
of linear operators on $\IR^n$ which contains $-I$  can be realized as the 
group of isometries of some sequence space (see \cite{GL79}), we will not 
attempt here to characterize them completely.    

We prove that if $X, \ Y$ are symmetric finite dimensional sequence  spaces, 
i.e., spaces with symmetric norms, then isometries of $X(Y)$ necessarily 
preserve the cross product structure except when $X = \ell_p$ and $Y$ can  be decomposed as an $\ell_p$-direct sum of two nonzero  
subspaces.  All other 
possible isometries in the exceptional case are also characterized. 
It is worth noting that even in this special type of cross product spaces,
the results in the complex case and the real case are quite different
when $\dim Y = 2$ or $4$ (cf. Corollary 3.4 and Theorem 4.1).

\smallskip
For simplicity of notation, we shall always assume that we have a
normalized absolute norm, i.e., all standard unit vectors
have norm 1.

\bigskip\noindent
{\bf 2. Intrinsic Cross Product Structure}\rm

We begin with the definition of fibers which is modelled on the  structure 
of the space $X(Y_1,Y_2,\ldots,Y_k)$, where each of the $Y_i$ is a ``fiber''. 

\smallskip\noindent
{\bf Definition 2.1} \
Let $X$ be a sequence space with a normalized absolute norm $N$, and  
let $\{e_1, \dots, e_k\}$ be the corresponding 1-unconditional basis,  
where $k\in\IN \cup \{\infty\}$.  A non-empty proper subset $S$ of 
$\{1,\ldots,k\}$ is called a {\it fiber} if for all $a_s,a'_s\in  
\IF$, $s\in S$ 
$$N\biggl( \sum_{s\in S} a_s e_s\biggr) 
= N\biggl( \sum_{s\in S} a'_s e_s\biggr)$$
implies that for all $b_i\in \bbF$, $i\in \{1,\ldots,k\}\setminus S$ 
$$N\biggl( \sum_{s\in S} a_s e_s + \sum_{i\notin S} b_i e_i\biggr) 
= N\biggl( \sum_{s\in S} a'_s e_s + \sum_{i\notin S} b_i e_i\biggr)\  
.$$
Moreover, the corresponding {\it fiber space} is defined by
$$X_S = \span\{e_s: s \in S\}.$$

\medskip
Here we mention a few examples of fibers.
\begin{enumerate}
\item Clearly, in any $X$, if $S$ is a singleton then $S$ is a fiber. 

\item In $\ell_p$, $1\le p\le \infty$, every  non-empty proper subset
$S$ of $\{1,\ldots,k\}$ is a fiber. 

\item  Let $1 \le p \le \infty$.  Let $\ell_p^k$ be the  
$k$-dimensional $\ell_p$ space, and let $\ell_p^k(Y)$, where
$Y$ cannot be written as a direct sum of $Y_1, \dots, Y_r$ such 
that $Y=\ell_p(Y_1, \dots, Y_r)$ (i.e. $Y$ has no nontrivial  
$\ell_p-$summands).   Then fiber spaces of $\ell_p^k(Y)$ are of
the form $W_1 \times \cdots \times W_k$  where $W_i$ equal $\{0\}$ or  
$Y$  for all $i= 1,\dots, k,$ or $W_i=\{0\}$ for all $i$ except exactly  
one, say $i_0$, and $W_{i_0}$ is a fiber space in $Y$.

\item In the above example, if the space $\ell_p^k$ is replaced by a 
different space $X$, as shown below,  we do not need any assumptions 
on $Y$ (we  even allow $Y=X$). We will see that  $W_1 \times \cdots \times 
W_k$ is a fiber  space in $X(Y)$ if and only if $W_i=\{0\}$ for all $i$ 
except exactly one, say  $i_0$, and $W_{i_0}$ is a fiber space in $Y$. 
\end{enumerate}

\smallskip
We want to analyze fibers on $X$ which are maximal with respect to 
inclusion. We start with the following observation: 

\smallskip\noindent
{\bf Proposition 2.2} \ \sl Let $X$ be a $k$-dimensional sequence  
space with a normalized absolute norm $N$.
Suppose there exist two maximal fibers $S,T$ such that $S\cap T\ne
\emptyset$. Then $S\cup T = \{1,\ldots,k\}$ and
$X= \ell_p (X_{T\setminus S},X_{T\cap S},X_{S\setminus T})$. 

Proof.
\rm 
Suppose $S$ and $T$ are maximal so $S\setminus T\ne \emptyset$ and 
$T\setminus S\ne\emptyset$. Assume that $i_0\in S\cap T$. 
Since $S$ is a fiber 
\begin{equation}\label{non1}
N\biggl( \sum_{s\in S\cup T} a_s e_s\biggr) 
= N\left( N\biggl( \sum_{s\in S} a_s e_s\biggr) e_{i_0} 
+ \sum_{s\in T\setminus S} a_s e_s\right) .
\end{equation}
Moreover,
\begin{eqnarray*}
& ~~~ &
N\biggl( \sum_{s\in S\cup T} a_s e_s + \sum_{i\notin S\cup T} b_i  
e_i\biggr)  \cr
& = &  N\left( N \biggl( \sum_{s\in S} a_s e_s \biggr) e_{i_0} 
+ \sum_{s\in T\setminus  S} a_s e_s 
+ \sum_{i\notin S\cup T} b_i e_i\right) \qquad\qquad\qquad 
  \mbox{since $S$ is a fiber} \cr
& = & N\left( N\left( N\biggl( \sum_{s\in S} a_s e_s\biggr) e_{i_0} 
+ \sum_{s\in T\setminus S} a_s e_s\right)  e_{i_0} 
+ \sum_{i\notin S\cup T} b_i e_i\right) \qquad 
\mbox{since $T$ is a  fiber} \cr
& = & N\left( N \biggl( \sum_{s\in S\cup T} a_s e_s\biggr) e_{i_0} 
+ \sum_{i\notin S\cup T} b_i e_i\right)\ .  \hskip 1.6in
\mbox{by (\ref{non1})} \cr
\end{eqnarray*}
Therefore $S\cup T$ is a fiber, and by maximality of $S$, $S\cup T= 
\{1,\ldots,k\}$. Take indices $s_0\in S \setminus T$, $i_0\in S\cap T$, 
$t_0\in  T\setminus S$. Using consecutively the fact that $T$ and $S$ are 
fibers we get for  all scalars $x_1,x_2$: 
\begin{equation}\label{non2}
N(x_1e_{s_0} + x_2 e_{i_0}) 
= N(x_1 e_{s_0} + x_2 e_{t_0}) 
= N(x_1 e_{i_0} + x_2 e_{t_0})
\end{equation} 
Next, since $S$ is a fiber and using (\ref{non2}) we get 
\begin{eqnarray}
N(x_1 e_{s_0} + x_2 e_{i_0} + x_3 e_{t_0}) 
& = & N (N (x_1 e_{s_0} + x_2 e_{i_0}) e_{s_0} + x_3 e_{t_0} )  \cr
& = & N (N (x_1 e_{s_0} + x_2 e_{i_0}) e_{s_0} + x_3 e_{i_0})  
\label{non3} 
\end{eqnarray}
Similarly, since $T$ is a fiber and by (\ref{non2}) 
\begin{eqnarray}
N (x_1 e_{s_0} + x_2 e_{i_0} + x_3 e_{t_0}) 
& = & N(x_1 e_{s_0} + N(x_2 e_{i_0} + x_3 e_{t_0}) e_{i_0})  \cr
& = & N(x_1 e_{s_0} + N(x_2 e_{s_0} + x_3 e_{i_0}) e_{i_0})  
\label{non4} 
\end{eqnarray}
Now let $f: \bbR^2 \lra \bbR$ be defined by $f(x_1,x_2) = 
N(x_1 e_{s_0} + x_2 e_{i_0})$. Then (\ref{non3}) and (\ref{non4})  
will take the form: 
$$f(f(x_1,x_2),x_3) = f(x_1 ,f (x_2,x_3))\ .$$
By a theorem of Bohnenblust \cite{Bo}, there exists $p$ with $1\le  
p\le \infty$  such that
$$
f(x_1,x_2) = \cases{
(|x_1|^p + |x_2|^p)^{1/p} & if $p<\infty$ \cr
\max (|x_1| ,|x_2|) & if $p=\infty$\ .}$$
Therefore for all scalars $x_1,x_2,x_3$ 
\begin{equation} \label{non5}
N(x_1 e_{s_0} + x_2 e_{i_0} + x_3 e_{t_0}) = \ell_p (x_1,x_2,x_3)\ .
\end{equation} 
By the fact that $S$ and $T$ are fibers and by (\ref{non5}),  
for any $\{x_s\}_{s=1}^k$, we have:
\begin{eqnarray*}
N\biggl( \sum_{s\in T} x_s e_s\biggr) 
& = & N\left( N\biggl( \sum_{s\in S\cap T} x_s e_s\biggr) e_{s_0} 
 + \sum_{s\in T\setminus S} x_s e_s\right)  \cr
& = & N\left( N\biggl( \sum_{s\in S\cap T} x_s e_s\biggr) e_{s_0} 
 + N\biggl( \sum_{s\in T\setminus S} x_s e_s\biggr) e_{t_0}\right)  
\cr
& = & \ell_p \left( N\biggl( \sum_{s\in S\cap T} x_s e_s\biggr) , 
 N\biggl( \sum_{s\in T\setminus S} x_s e_s\biggr)\right)\ . \cr
\end{eqnarray*}
Similarly: 
\begin{eqnarray*} 
N\biggl( \sum_{s=1}^k x_s e_s\biggr) 
& = & N\left( N\biggl( \sum_{s\in T} x_s e_s\biggr) e_{t_0} 
 + \sum_{s\in S\setminus T} e_s x_s\right)  \cr
& = & N\left( N\biggl( \sum_{s\in T} x_s e_s\biggr) e_{t_0} 
 + N\biggl(\sum_{s\in S\setminus T} e_s x_s\biggr) e_{s_0}\right) \cr
& = & \ell_p \left( N\biggl( \sum_{s\in T} x_s e_s\biggr) , 
 N \biggl( \sum_{s\in S\setminus T} e_s x_s\biggr) \right)  \cr
& = & \ell_p \left( N\biggl( \sum_{s\in T\setminus S} x_s e_s\biggr) ,
 N \biggl( \sum_{s\in T\cap S} x_s e_s\biggr) , 
 N \biggl( \sum_{s\in S\setminus T} e_s x_s\biggr) \right)\ . \cr
\end{eqnarray*}
Thus $X= \ell_p (X_{T\setminus S}, X_{T\cap S}, X_{S\setminus T})$.
\qed

Thus maximal fibers determine the structure of a sequence space $X.$  
To prove our main theorem on the cross product structure, we need a 
definition of a special 2-dimensional real space different form $\ell_p^2,$ 
which can be decomposed into $\ell_p$ sum of its nonzero subspaces  
(see \cite{LW}): 

\smallskip\noindent
{\bf Definition 2.3} \
Let $1\le p<\infty, \ p\neq2.$ Let $E_p(2)$ denote the space $\bbR^2$  
with the following norm: 
$$\|(x,y)\|_{E_p} = \left( {|x+y|^p\over 2} +  
{|x-y|^p\over2}\right)^{1/p}.$$
If $p=\infty$ define:
$$\|(x,y)\|_{E_{\infty}} = \max( {|x+y|},{|x-y|}) =  
\|(x,y)\|_{\ell_1}.$$

\smallskip
Observe that $E_p(2)$ is isometric to $\ell_p^2$ through the isometry 
$T:E_p(2)\lra \ell_p^2$ defined by $T(x,y) = 2^{-1/p}(x+y,x-y).$ 

\smallskip\noindent
{\bf Theorem 2.4} \sl
Let $X$ be a sequence space with a normalized absolute norm $N$, and  
let $\{e_1, \dots, e_k\}$ be the corresponding 1-unconditional basis,  
where $k\in\IN \cup \{\infty\}$. Then there exists a partition $S_1, \dots, 
S_m$ of $\{1, \dots, k\}$ such that $X$ is a direct sum of 
$X_i = \span\{ e_s: s \in S_i\}$, $1 \le i \le  m$, and one of the 
following holds. 
\begin{enumerate}
\item[{\rm (i)}]
Every $S_i$ is a maximal fiber and $X = X_0(X_1, \dots, X_m)$
where the norm $N_0$ on $X_0$ is defined by
$$N_0(a_1,\ldots,a_m) = N\left(\sum_{i=1}^m a_i e_{s_i}\right)$$ 
for some $s_i\in S_i$. 

\item[{\rm (ii)}]
There exists $p$ with $1\le p\le\infty$ such that
$X = \ell_p(X_1,\dots,X_m),$ where 
\begin{itemize}
\item[{\rm (a)}]
at most one of the $X_i$'s equals $\ell_p^{\dim X_i}$,

\item[{\rm (b)}]
some of the $X_i$'s equal to $E_p(2)$ if $\IF = \IR$ and $p \ne 2$,

\item[{\rm (c)}]
and the rest of the $X_i$'s are such that  $\dim X_i \ge 2$ and
$X_i$ is not an $\ell_p$ sum of two nonzero subspaces.
\end{itemize}
\end{enumerate}

Proof.   \rm Suppose (i) does not hold. By Proposition 2.2, there  
exists $p$, $1\le p\le\infty$ and spaces $Y_1,Y_2,Y_3$ so that $X =  
\ell_p(Y_1,Y_2,Y_3)$ and $Y_i =\span\{e_s: s\in A_i\},$  \ for some 
partition  $\{A_i\}_{i=1}^3$ of $\{1,\dots,k\}.$   Among all the  
decompositions of the space $X$ into $\ell_p$ sum, let $R_1 \cup \cdots 
\cup R_s$ be a maximal partition  of $\{1, \dots, k\}$ so that 
$X = \ell_p(Z_1, \dots, Z_s)$ with $Z_i =  \span\{e_{r}: r \in R_i\}$. 
If $X$ is real and $p\neq 2$, then  for each $1 \le i \le s$ we 
have one of the three possibilities (cf. \cite{LW}): 
\begin{itemize} 
\item[(a)] $R_i$ is a singleton.

\item[(b)] $R_i$ has two elements and $Z_i = E_p(2)$.

\item[(c)] $R_i$ has at least two elements and $Z_i$
cannot be decomposed as an $\ell_p$-direct sum of two nonzero  
subspaces.
\end{itemize}
If $X$ is complex or if $p=2$ only (a) and (c) can happen (cf.  
\cite{BL77}).

Let $S_1$ be the union of the $R_i$ which are singletons if they  
exist, and rename the other $R_i$ as $S_j$ if necessary.  We see that 
condition (ii) holds.
\qed

\medskip
Several remarks are in order in connection with Theorem 2.4.
\begin{enumerate}
\item
In both cases (i) and (ii) we present $X$ as a cross product (or  
direct sum) of disjoint factors (or factor subspaces) $X_1,\ldots,X_m$.
Notice that $X_1,\ldots, X_m$ are uniquely determined by $X$ and that 
each of the spaces $X_i$ may be  further decomposable into 
factors (we do not consider factors of $X_i$'s as factors of $X$, 
sometimes we will call them second generation factors of $X$). 

\item
Notice also that in case (i), the factors are precisely maximal  fibers of 
$X$. In case (ii) factors are uniquely determined by maximal fibers.  Namely, 
$X_1 = \ell_p$, then it  has supports on the union of those points $\{i\}$ in 
$\{1,\ldots,k\}$ such that $\{1,\ldots,k\}\setminus \{i\}$ is a maximal fiber; 
and each of the other $X_j$ with $\card (\supp X_j)\ge2$ has
supports on a complement of a maximal fiber. 

\item
Evidently, if $X$ has an explicit cross product structure i.e. 
if $X=Y(Y_1,\dots,Y_m)$ then Theorem 2.4 can be used to
regroup the factors of $X$ so that conditions (i) or (ii) of
Theorem 2.4 holds.  If no regrouping is necessary, we say that $X$  
has {\it reduced cross product structure}.
\end{enumerate}

\bigskip\noindent
{\bf 3. Isometries of Complex Sequence Spaces}

Before further analysis of isometries of $X$ we need to introduce 
another definition. After \cite{ST} (cf. also \cite{KW}) we define 
an equivalence relation $\sim$ on the indices $\{1,\ldots,k\}$. 
We say that $s\sim t$ if $N(\sum_{i=1}^k a_i e_i)=N(\sum_{i=1}^k b_i e_i)$ 
whenever $\ell_2 (a_s ,a_t) = \ell_2 (b_s,b_t)$ and $a_i = b_i$ for all 
$i\ne s,t$. 

Schneider and Turner showed that $\sim $ is indeed an equivalence  
relation and that equivalence classes of $\sim$ are isometrically isomorphic 
to $\ell_2$ (with appropriate dimension) \cite[Lemma 2.3]{ST}. 
Clearly equivalence classes of $\sim$ are fibers in $X$, and hence  
they are contained in maximal fibers of $X$. 
We will call equivalence classes of $\sim$ maximal $\ell_2$-fibers. 
Notice that every subset of a maximal $\ell_2$-fiber is a subfiber  
and it is a (non-maximal) $\ell_2$-fiber. 

The results in \cite{ST} and \cite{KW}
state that every isometry of $X$ preserves maximal 
$\ell_2$-fibers. This fact has very important consequences for us. 
Namely we have: 

\smallskip\noindent
\bf Proposition 3.1
\sl 
Let $X$ be a  complex sequence space with  a normalized absolute basis
$\{e_1, \dots, e_k\}$, where $k\in\IN \cup \{\infty\}$. 
If $T$ is a surjective isometry of $X$, then $T$ 
preserves maximal fibers of $X$, i.e.,  if 
$S\subset \{1,\ldots,k\}$ is a maximal fiber of $X$ 
then $\supp(T(\span \{e_i:i\in S\}))$ is a maximal fiber of $X$. 

Proof.
\rm
Denote by $\{J_j\}_{j\in \Lambda}$ $(m\le\infty)$  the collection of 
all maximal $\ell_2$-fibers in $X$. 
Then $\{1,\ldots,k\}= \bigcup_{j\in\Lambda} J_j$. 
By \cite[Theorem 6.1]{KW} there exists a permutation $\sigma$ of 
$\{1,\ldots,k\}$ such that for all $j\in\lambda$ 
\begin{equation} \label{KW}
\supp (T(\span \{ e_s : s\in J_j\})) = J_{\sigma (j)}.
\end{equation}
Thus we can define a map $\w T$ which operates on maximal  
$\ell_2$-fibers and their unions by 
$$\w T (J) = \supp (T (\span  \{e_s : s\in J\})) =  
\bigcup_{j\in\Lambda_1} J_{\sigma (j)}$$ 
where $J= \bigcup_{j\in\Lambda_1} J_j$ for any $\Lambda_1\subset  
\Lambda$. 

Let $S\subset \{1,\ldots,k\}$ be a fiber in $X$. 
Then $S= \bigcup_{j\in \Theta} J_j$ for some $\Theta \subset\Lambda$. 
Thus $\w T (S)$ and $\w T (S^c)$ are well defined and disjoint. 

Now let $a,a'\in X$ be such that $\supp a \cup \supp a'\subset \w  T(S)$ 
and $N(a) = N(a')$, and let $b\in X$ with $\supp b\subset \w T(S^c)$. 
Then $\supp T^{-1} (a) \cup \supp T^{-1}(a') \subset S$ and 
$\supp T^{-1}(b) \subset S^c$. Thus, since $S$ is a fiber, we have 
$$
N(a+b) =  N(T^{-1} (a) + T^{-1} (b)) = N(T^{-1} (a') + T^{-1}(b)) 
=  N(a'+ b)\ . $$
Therefore $\w T (S)$ is a fiber in $X$. 
Moreover, this implies that $({\w T}^{-1})$ can be applied to the set 
$\w T(S)$ and $({\w T}^{-1})(\w T(S)) = S$. 

Now if $S$ is a maximal fiber then  $\w T(S)$ is also a maximal fiber. 
Indeed, if $\w T (S)$ is not maximal, say $\w T(S)$ is a subfiber of a 
proper fiber $S_1$ then $({\w T}^{-1})(S_1)$ is a proper fiber in $X$  
which contains $({\w T}^{-1})(\w T (S))=S$ which contradicts the
maximality of $S$. 
\qed

We are now ready to present our main result about the form of isometries 
of complex sequence spaces with 1-unconditional basis. 

\smallskip\noindent
\bf Theorem 3.2 \sl
Let $X$ be a complex sequence space with 1-unconditional basis 
$\{e_1,\ldots,e_k\}$, 
$k\in\IN\cup\{\infty\}$, and let $X=X_0(X_1,\ldots,X_m)$ where  
$X_1,\ldots,X_m$ are factors as described in Theorem 2.4.
Then $T$ is a surjective isometry of $X$ if and only if there exists a 
permutation $\pi$ of $\{1,\ldots,m\}$ such that the norm $N_0$ on $X_0$ 
satisfies $N_0(z_1,\ldots,z_m) = N_0(z_{\pi (1)}, \ldots, z_{\pi (m)})$ 
and there exists a family of surjective isometries $S_j: X_{\pi (j)}  
\lra X_j$ such that 
\begin{equation}\label{basic}
T(x_1,\ldots, x_m) = (S_1x_{\pi (1)} ,\ldots, S_m x_{\pi (m)}) 
\end{equation} 
for all $(x_1,\ldots,x_m) \in X_0(X_1,\ldots, X_m)=X$.

Proof. \rm
The result follows quickly from Proposition~3.1 and Remark~2 after 
Theorem~2.4. The only thing that needs to be verified is that if $S$ is a 
maximal  fiber such that $\card (S^c)=1$ then $\card (\w T (S^c))=1$. But this 
is clear since then $S^c$ is a maximal $\ell_2$-fiber of cardinality~1 and by 
\cite[Theorem 6.1]{KW} $\w T(S^c)$ is a maximal $\ell_2$-fiber of 
cardinality~1. Thus $\card (\w T (S)^c) = \card (\w T(S^c))=1$. Therefore, if 
$X$ has structure described in Theorem 2.4 (ii), then $X_1$ is mapped to 
$X_1$, i.e., $\pi (1)=1$. \qed 

We would like to make two remarks:
\begin{itemize}
\item[1.]
Proposition~3.1 is also valid for a real sequence space, whose isometry group 
is contained in the group of signed permutations.  Such a space has no 
nontrivial $\ell_2-$fibers and (\ref{KW}) holds, so the conclusion of 
Proposition~3.1 will follow. 

\item[2.]
Theorem 3.2 is valid in those real spaces for which  Proposition~3.1  holds.  
In particular, by the discussion in the preceding paragraph, Theorem 3.2 is 
valid for a real sequence space, whose isometry group is contained in the 
group of signed permutations. 
\end{itemize} 

Theorem 3.2 provides a complete description of surjective 
isometries of complex sequence spaces. Below we present some immediate 
corollaries about the form of isometries of spaces with explicit cross
product structure (cf. Remark 3 after Theorem 2.4). 

\medskip\noindent
{\bf Corollary 3.3} \sl
Let $X=Z(X_1,\ldots,X_m)$ be the space with explicit reduced cross product 
structure. Then every surjective isometry of $X$ onto itself has form 
(\ref{basic}).  \rm

\medskip\noindent
{\bf Corollary 3.4} \sl
%\begin{cor}\label{xy}
Let $X,Y$ be complex sequence spaces not both equal to $\ell_p$ (with 
the same $p$). Then the isometries of $X(Y)$ have form (\ref{basic}).  
\rm

\medskip
Notice that if factors $X_1,X_2,\ldots,X_m$ can be further decomposed into 
simpler second generation factors (as mentioned in Remark 1 after Theorem 
2.4), then one can again use Theorem~3.2 to conclude that isometries 
$S_1,\ldots,S_m$ have form  $(8)$. In particular, one can inductively 
describe isometries of spaces of the form 
$X_1(X_2(\ldots(X_m)\ldots)) $, where $X_1,X_2,\ldots,X_m$ are
complex sequence spaces such that for any $i=1,\ldots,m-1,$ $X_i$ and 
$X_{i+1}$ are not simultaneously equal to $\ell_p$ with the same $p$. 
We leave the exact statement to the interested reader.
It is interesting to note that the group of isometries of   
$X_1(X_2(\ldots(X_m)\ldots))$
does not depend on entire isometry groups of  $X_1,\ldots,X_{m-1}$,  
but only on intersection of these groups with the group of signed permutation  
operators and the isometry group of $X_m$.

\newpage
\bigskip\noindent
{\bf 4. Isometries of Real Sequence Spaces} 

The description of isometries of real sequence spaces is more  complicated 
than in the complex case. The main difference is in classification of spaces 
whose group of isometries  is  contained in the group of signed  permutations. 
In the complex case Schneider, Turner \cite{ST} and Kalton, Wood  \cite{KW} 
determined that group of isometries  is contained in the group of signed 
permutations if and  only if the space does not have nontrivial 
$\ell_2$--fibers. 

In the real case similar classification is not valid. In fact, we have the 
following examples of spaces which do not contain any copies of $\ell_2$ and 
which allow non-disjointness preserving isometries, i.e., isometries  
that are not signed permutation operators.  
\begin{itemize}
\item[1.]
%\begin{example}\label{exep}
Let $p\neq 2$ and let $E_p(2)$  be the space defined in Definition 2.3.
Then $E_p(2)$ is isometric to $\ell_p^2$ through the isometry 
$T:E_p(2)\lra\ell_p^2, \ $ $T(x,y) = 2^{-1/p}( {x+y } , {x-y}).$ Let  
$X= Y(\ell_p^2,E_p)$ for any 2-dimensional real symmetric space $Y$, and  
let $S:X\lra X$ be an isometry defined by $S(v_1,v_2) = (Tv_2 ,T^{-1} v_1)$ 
where $v_1 \in \ell_p$, $v_2\in E_p$. Then $S$ is a non-disjointness  
preserving isometry of $X$. 
%\end{example} 

\item[2.]
%\begin{example} \label{exex}
Let $X$ be any 2-dimensional real symmetric space, $X\neq\ell_2.$ Put 
$$\|(x,y)\|_{E_X} = {1\over \|(1,1)\|_X} \|(x+y,x-y)\|_X\ .$$ 
Consider $Z= Y(X, E_X)$, where $Y$ is any 2-dimensional real symmetric space. 
Then, similarly to Example 1, there exists a non-disjointness preserving 
isometry of $Z$. 
%\end{example} 

\item[3.]
%\begin{example} \label{exmultex}
Spaces in Example 2 can be generalized to higher dimensions by 
taking any spaces $X_1, X_2$ which are isometric through a non-disjointness
preserving isometry  (spaces like that can be constructed e.g. by taking 
direct products of $X$ and $E_X$, cf. also \cite[Theorem~4]{Rseq}).  
Then let 
$Z= Y(X_1,  X_2)$ for any  symmetric space $Y$. 
%\end{example} 
\end{itemize}

\smallskip
The above examples show  isometries which are not 
signed permutations but which nevertheless ``preserve cross product 
structure'', i.e., have canonical form (\ref{basic}). One would hope  
that this is always true, however the following examples show 
the contrary. 
\begin{itemize}
\item[4.]
%\begin{example}\label{exlpep}
Consider $\ell_p(\ell_p^3,E_p(2)) = (\IR^3\times \IR^2, N)$ with
$$N(x_1,x_2,x_3,y_1,y_2) = \ell_p(\ell_p(x_1,x_2,x_3),E_p(y_1,y_2))\}.$$  
As described in Example~1, $\ell_p^2$ is isometric with $E_p(2)$
via the isometry $T$.  
Thus one can define isometry $S$ on $\ell_p(\ell_p^3,E_p(2))$ by
$$
S(x_1,x_2,x_3,y_1,y_2)
= \left(T(y_1,y_2), x_3, T^{-1}(x_1, x_2)\right).$$
Clearly $S$ does not have form $(8)$.
%\end{example}
\end{itemize}

\smallskip
Surprisingly, similar pathology is possible even in spaces of the  
form $X(Y)$, where $X, \ Y$ are symmetric.
\begin{itemize}
\item[5.]
Consider the space $\ell_p^2(E_p(2)) = \ell_{p}(E_p(2),E_p(2))$. Then we have 
$$\matrix{
& ~~~  N(x_1,x_2,y_1,y_2) \hfill \cr  & \cr
&= \ell_p(E_p(x_1,x_2),E_p(y_1,y_2)) \hfill \cr & \cr
&= \ell_p(2^{-1/p} \ell_p(x_1+x_2, x_1 - x_2),
     2^{-1/p} \ell_p(y_1+y_2,y_1 - y_2)) \hfill \cr  & \cr
&= 2^{-1/p} \ell_p(x_1+x_2, x_1 - x_2, y_1+y_2,y_1 - y_2) \hfill\cr & \cr
&= 2^{-1/p} \ell_p(\ell_p(x_1+x_2, y_1+y_2), 
    \ell_p(x_1 - x_2,y_1 - y_2)) \hfill\cr & \cr
&= \ell_p\biggl(E_p\biggl(2^{-1}(x_1+x_2+y_1+y_2,x_1+x_2-y_1-y_2)\biggr), 
\hfill\cr
& ~~~~~ E_p\biggl(2^{-1}(x_1 - x_2+y_1 - y_2,x_1 - x_2-y_1+y_2)\biggr)\biggr) 
\hfill\cr
&=  
N\biggl(2^{-1}(x_1+x_2+y_1+y_2,x_1+x_2-y_1-y_2,  \hfill\cr
& ~~~~  x_1 - x_2+y_1 - y_2,x_1 - x_2-y_1 +y_2) \biggr) \cr
}$$
Thus the linear map defined by 
\begin{eqnarray*}
S(1,0,0,0) &=& 2^{-1}(1,1,1,1) 
\quad \quad \qquad S(0,1,0,0) =  2^{-1}(1,1,-1,-1) \cr
S(0,0,1,0) &=& 2^{-1}(1,-1,1,-1) 
\qquad \ S(0,0,0,1) =  2^{-1}(1,-1,-1,1)
\end{eqnarray*}
is an isometry, and clearly $S$ does not preserve disjointness of  
vectors. This isometry in the case when $p=\infty,$ i.e.
$X= \ell_\infty^2(\ell_1^2)$, and $p=1$, i.e. $X= \ell_1^2(\ell_\infty^2)$ is 
described in \cite[Theorem~3.2(b)]{KL}. Notice that $E_p(2)$ can  be decomposed as an $\ell_p$-direct sum of two nonzero  
subspaces. Thus this example is consistent with \cite[Proposition~2]{Gr}, which sais that if $E$ is an $\ell_p-$sum of two nonzero subspaces then there exists an isometry of $\ell_p(E)$ which is not of the form $(9)$.
 It is interesting that this is, 
in fact, the only possible example as shown in 
Theorem 4.1 below.
%\end{example}
\end{itemize}

\smallskip
It becomes of interest to characterize spaces with cross product structure, 
which is preserved under action of all isometries. 

First we list classes of spaces whose isometry group is contained in 
the group of signed permutations and thus Theorem 3.2 can be applied to 
conclude that indeed the cross product structure  is preserved by all 
isometries. This holds for: 
\begin{itemize}
\item[(1)] spaces with $\Delta$-bases (\cite{GL79})

\item[(1a)] In particular for spaces of the form  $Z(X_1,X_2,\ldots,X_k)$ 
where $k\le\infty$, $\dim X_i \in \bbN \setminus\{ 2,4\}$ and $X_i$ is 
a symmetric space not equal to $\ell_2$, $i=1,\ldots,k$, $Z$ is  
arbitrary. 

\item[(2)] spaces which are $p$-convex with constant 1 for  
$2<p<\infty$ (\cite{Rseq}),

\item[(2a)]  spaces which are strictly monotone, smooth at every basis 
vector and $q$-concave with constant 1 for $1<q<2$ (\cite{Rseq}). 
\end{itemize}

\noindent
Further Rosenthal \cite{R86} (cf. also \cite{Rxh}) showed that 
Theorem 3.2 holds in real spaces of the form: 
\begin{itemize}
\item[(3)] $Z(X_1,\ldots,X_k)$ where $k\le \infty$, $\dim X_i\ge 2$  
and $X_i=\ell_2$ for all $i=1,\ldots,k$, $Z\ne\ell_2$.
\end{itemize}

\smallskip
Below we study the group of isometries of spaces of the form $X(Y)$ 
where $X,Y$ are finite-dimensional symmetric spaces with  $\dim Y\ge 2$. 
This class of spaces has a sizable intersection with class (1a) and  
(3), but we do allow $\dim Y$ to be 2 or 4, which are excluded by (1a). 
Thus  we allow the situation when the isometry group  is not contained in the 
group of signed permutations. Our proof is somewhat technical but it is more 
elementary compare with the one in \cite{GL79}.  It is possible that our
different approach may lead to some insight to the general problem.
So we present the entire proof including the previously proven cases.

\medskip\noindent
{\bf Theorem 4.1} \sl Suppose $\dim X = n$ and $\dim Y = m$,
and not both $N_1$ and $N_2$ are $\ell_p$ $($with the same $p)$.
If $\Psi$ is an isometry for $N$ on $X(Y)$, then 
\begin{itemize}
\item[{\rm (i)}] $\Psi$ is of the form 
\begin{equation}\label{new}
(y_1, \dots, y_n) \mapsto (S_1 y_{\pi(1)}, \dots, S_k y_{\pi(n)})
\end{equation}
for some isometries $S_i$ of $Y$ and some permutation $\pi$ of 
$\{1, \dots, n\}$, or 

\item[{\rm (ii)}] $X = \ell_p^n$,   $Y = E_p(2)$ for some $p, \ \ 1\le p\le\infty, \ \  p\neq2,$ 
and  
\begin{equation}\label{new2}
 \Psi  {\rm\ \  permutes\ \   the\ \   matrices\ \   in \ \  the \ \   
set\ \  }
\{\pm (e_1^Y \pm e_2^Y) (e_j^X)^t: 1 \le j \le n\},
\end{equation}
where $e_i^Y, \ e_i^X$ denote column basis vectors in $Y$ and $X$, resp.
\end{itemize}
\rm

\medskip
Notice that, clearly, if $N_1 = N_2 = \ell_p$ then $X(Y) = \ell_p^{mn}$ and 
the isometry group is well known. 
  
The choices of $S_i$ in Theorem 4.1 (i)
are very restrictive. If $N_2 = \ell_2$, the $S_i$ is orthogonal (on $Y$); 
otherwise, $S_i$ a signed permutation operator on $Y$ unless for $m = 4$ 
or $m = 2$.  When $m = 4$, there are a few more possibilities for $S_i$.  
One may, for example,  see \cite{Rol} (see also \cite{DLR,Br}) for more 
details. This exceptional case will be treated separately in our proof.  When 
$m = 2$, $S_i$ must be chosen from a dihedral group. 

In the following discussion, we shall 
identify $X(Y)$ with the space $\IR^{m\times n}$
of $m\times n$ real matrices, and identify the linear operator
$\Psi$ on $X(Y)$ with its matrix representation 
relative to the standard basis
$$\{e_{11}, e_{21}, \dots, e_{m1}, \dots, e_{1n}, e_{2n}, \dots,  
e_{mn}\}.$$
We shall also use the following notations in our discussion.

$O(m)$: the orthogonal group on $\IR^m$.

$O(mn)$: the orthogonal group on $\IR^{m\times n}$.

$P(m)$: the group of $m\times m$ permutation matrices.

$GP(m)$: the group of $m\times m$ signed permutation matrices.

$GP(mn)$: the group of linear operators on $\IR^{m\times n}$ that
permute and change the signs of the entries of $A \in \IR^{m\times n}$.

$P\otimes Q$: the tensor product of the matrices $P$ and $Q$ given
by $(P_{ij}Q)$.

$e_{ij}^{(n)}$: the $n\times n$ matrix with one at the $(i,j)$  
position and zero elsewhere.

$E_{pq}^{(ij)}$: the $mn\times mn$ matrix $e_{ij}^{(n)} \otimes 
e_{pq}^{(m)}$.

We shall also use the concept of the {\it Wreath product} of two  
groups of linear operators.  For simplicity, we consider the special case
when $G$ is a group of linear operators (identified as matrices)
acting on $\IR^m$.
The Wreath product of $G$ and $P(n)$, denoted by $G*P(n)$, is the 
group of linear operators on $\IR^{m\times n}$  of the form
$$ [y_1| \cdots| y_n] \mapsto [U_1(y_1)| \cdots | U_n(y_n)]V $$
for some $U_1, \dots, U_n \in G$ and $V \in P(n)$.

With this definition, Theorem 4.1 implies that the
isometry group of $N$ is the Wreath product of the isometry group
of $N_2$ and $P(n)$ if $m \ne 2$. In particular, isometries will 
always preserve the cross product structure of $X(Y)$.

\medskip
It is also interesting to note that in our proof, we actually  determine
all possible closed overgroups $G$ of $H = GP(m)*GP(n)$ in $O(mn)$. In 
particular, if $G$ is infinite then $G = O(m)*GP(n)$ or $O(mn)$; 
if $G$ is finite then $G$ is one of the following:
\begin{itemize}
\item[(a)] $GP(m)*GP(n)$: the Wreath product of $GP(m)$ and $GP(n)$,

\item[(b)] $GP(mn)$: the group of signed permutations on $\IR^{m\times n}$,

\item[(c)] $A_{2n}$: the Weyl group of $A_{2n}$ type realized as
the group of orthogonal operators on $\IR^{2\times m}$ 
that permute the set $\{\pm(e_1^Y\pm e_2^Y) (e_j^X)^t : 1 \le j \le n\}$
if $m = 2$, 

\item[(d)] $D_h*GP(n)$: the Wreath product of the dihedral group $D_h$ and  
$GP(n)$ if $m = 2$,

\item[(e)] $F_4$ and the normalizer of $F_4$: the Weyl group  of $F_4$, if 
$m =  n = 2$. 
\end{itemize}

\smallskip\noindent
In (e), there are two possible realizations of $F_4$, namely,
an overgroup of $GP(4)$ (e.g., see \cite{DLR}) or
the group generated by $H$ and $L_1$ mentioned in the proof of Lemma 4.7.
However, only the first realization can be an isometry group of X(Y).

\medskip
Our proof of Theorem 4.1 uses the basic ideas in \cite{DLR} 
(cf. also \cite{Br}) and some  
intricate arguments.  It would be nice to have a shorter conceptual  
proof. We begin our proof with the following corollary of Auerbach's 
Theorem \cite[Theorem~IX.5.1]{Rol} (see e.g. \cite[Theorem 2.3]{KL}).

\medskip\noindent
{\bf Lemma 4.2} \sl Let $G$ be the isometry group of $N$. Then
$G < O(mn)$, i.e., $G$ is a subgroup of $O(mn)$. \rm

\medskip 
We first deal with the case when the isometry group of $N$ is  
infinite.

\medskip\noindent
{\bf Lemma 4.3} \sl If the isometry group $G$ of $N$ is infinite,
then either $G = O(mn)$ or $G$ is the Wreath product $O(m)*P(n)$.

Proof. \rm By Lemma 4.2, $G$ is a subgroup of $O(mn)$. Since 
$G$ is closed and $O(mn)$ is a compact Lie group, $G$ is also a 
compact Lie group.  It is well-known that the Lie algebra of $O(mn)$
is $o(mn)$, the algebra of all skew-symmetric $mn \times mn$
matrices over real under the Lie product $[A,B] = AB-BA$.  
Suppose ${\bf g}$ is the Lie algebra of $G$.
Then ${\bf g}$ is a subalgebra of $o(mn)$.
Furthermore, by definition of $N$, we have $H:= GP(m)*GP(n) < G$
and $H$ acts on $G$ by conjugation, and so ${\bf g}$ is a
$H$-module under the action $(P,A) \mapsto P^tAP$ for any
$P \in H$ and $A \in {\bf g}$.  We shall show that there are only
two subalgebras of $o(mn)$ which are $H$-modules, and the two Lie
groups corresponding to the subalgebras are $O(mn)$ and $O(m)*P(n)$.

For any $A \in {\bf g}$ we write $A = (A^{(ij)})$ in $n\times n$  
block form
such that each block $A^{(ij)} \in \IR^{m\times m}$.  If there exists
$i < j$ such that $A^{(ij)} \ne 0$, we claim that ${\bf g} = o(mn)$.
First, note that there is $P \in P(n)$
such that the $(1,2)$ entry of $P^tZP$ is the $(i,j)$ entry of $Z$
for any $Z \in \IR^{n\times n}$.  Then $P\otimes I_m \in H$ and hence
$(P\otimes I_m)^t A (P\otimes I_m) \in {\bf g}$
will have nonzero $(1,2)$ block.  So, we may assume $(i,j) = (1,2)$.
Now suppose the $(p,q)$ entry of $A^{(12)}$ is nonzero.  Then
$D_1 = I_{mn} - E_{pp}^{(11)} \in H$ and hence
$A_1 = A - D_1AD_1 \in {\bf g}$.  Note that only the $p$th row and
$p$th column of $A_1$ are nonzero.  Now 
$D_2 = I_{mn} - E_{pp}^{(22)} \in H$ and hence
$A_2 = A_1 - D_2A_1D_2 \in {\bf g}$. Since $A_2$ is a nonzero 
multiple of $E_{pq}^{(12)} - E_{qp}^{(21)}$, it follows that
$Q^t(E_{pq}^{(12)} - E_{qp}^{(21)})Q \in {\bf g}$ for all $Q \in H$.  
As a result, ${\bf g}$ contains $E_{rs}^{(ij)} - E_{sr}^{(ji)}$
for all $1 \le i < j \le n$ and $1 \le s, r \le m$.
In particular, it contains the Lie product
$[E_{11}^{(12)} - E_{11}^{(21)}, E_{12}^{(21)} - E_{21}^{(12)}] =
E_{12}^{(11)} - E_{21}^{(11)}$.  It then follows that
$Q^t(E_{12}^{(11)} - E_{21}^{(11)})Q \in {\bf g}$ for all $Q \in H$.  
As a result, ${\bf g}$ also contains $E_{rs}^{(ii)} - E_{sr}^{(ii)}$
for all $1 \le i \le n$ and $1 \le r < s \le m$. So, ${\bf g}$  
contains
a basis of $o(mn)$ and we conclude that $G = O(mn)$.

Next, suppose all $A = (A^{(ij)}) \in {\bf g}$ satisfy 
$A^{(ij)} = 0$ if $i \ne j$.  Then we can show that
${\bf g}$ contains $E_{rs}^{(ii)} - E_{sr}^{(ii)}$
for all $1 \le i \le n$ and $1 \le r < s \le m$, by arguments
similar to the preceding case.  Taking the exponential map
for the elements of ${\bf g}$, we see that $G$ contains the
Wreath product $O(m)*\{I_n\}$.  Since $G$ also contains
$GP(m)*GP(n)$, we conclude that $O(m)*P(n) = O(m)*GP(n) < G$.
In particular, $N_2 = \ell_2$.  By the results in \cite{R86} 
(see also \cite{Rseq}),
we conclude that $G = O(m)*GP(n)$.
\qed

\medskip
Next we consider the case when the isometry group of $N$ is finite.  
We use the approach in \cite{DLR}, namely, determining all the finite 
overgroups of $GP(m)*GP(n)$ in $O(mn)$. We begin with the following  
lemma, which explains why one needs to exclude the cases when $n = 2, 4$ 
in \cite{GL79}.

\medskip\noindent
{\bf Lemma 4.4} \sl Suppose ${\rm dim}\,Y \ne 4, 2$.
If $G$ is the isometry group of $N$ and $G$ is finite, then
$G < GP(mn)$.

Proof. \rm For each $\Psi \in G$, we write $\Psi = (\Psi^{(ij)})$ in  
$n\times n$ block form such that $\Psi^{(ij)} \in \IR^{m\times m}$ for 
all $(i,j)$. Since $G < O(mn)$, $\Psi \in G$ implies $\Psi^t \in G$ as well. 

We use the technique in \cite{DLR} to show that for any $\Psi \in G$  
(in its matrix representation) the entries of $\Psi$ can only be $0, 1$ or  
$-1$. Suppose this is not true.  Let 
$$\mu = \min\{ a > 0: a \hbox{ is an entry of } \Psi \hbox{ for some } 
\Psi \in G\}.$$
Since $G < O(mn)$ by Lemma 4.2, we have $0 < \mu < 1$.  Denote by
$H = GP(m)*GP(n)$ as in the proof of Lemma 4.3.  Let $\Phi \in G$
have one of its entries equal to $\mu$.  Then there exist $P, Q \in  H$
such that the $(1,1)$ entry of $P\Phi Q$ equal $\mu$. Thus we may 
assume that the $(1,1)$ entry of $\Phi = (\Phi^{(ij)})$ equals $\mu$.
We consider several cases.

First, suppose $n$ is odd. By the
arguments in the proof of Theorem 2.4 in \cite{DLR}, there exists
$S \in GP(m)$ such that $(\Phi^{(11)})^tS\Phi^{(11)}$ has 
a positive $(1,1)$ entry equal to $\eta < \mu$.

Indeed, let $P = (P^{(ij)}) \in H$ be such that
$P^{(11)} = S$, $P^{(2k, 2k+1)} = I_m$ and 
$P^{(2k+1,2k)} = -I_m$ for $1 \le k \le (n-1)/2$, 
and $P^{(ij)} = 0$ for all other $(i,j)$.  Then
$\Phi^t P \Phi \in G$ has $(1,1)$ entry equal to $\eta = \mu^2 < \mu$,
which contradicts the definition of $\mu$.      

Next, suppose $n$ is even.  If $m$ is even, we can again obtain
$S \in GP(m)$ such that $(\Phi^{(11)})^tS\Phi^{(11)}$ has a positive
$(1,1)$ entry equal to $\eta < \mu$. Let  $P = (P^{(ij)}) \in H$
be such that $P^{(11)} = S$, 
$P^{(ii)} = \sum_{j=1}^{m/2} (e_{2j-1,2j}^{(m)} - e_{2j,  
2j-1}^{(m)})$
for $i = 2, \dots, n$, and $P^{(rs)} = 0$ for other $(r,s)$. Then
$\Phi^t P \Phi \in G$ has $(1,1)$ entry equal to $\eta < \mu$,
which contradicts the definition of $\mu$.

Finally, suppose $n$ is even and $m$ is odd.
Note that there exists $P \in H$ such that the first column
of $P\Phi$ is nonnegative.  So, we may assume that $\Phi$ has  
nonnegative column.  Furthermore, we may assume that
$\Phi_{11}^{(21)}$ attains the minimum among all entries
in the first columns of $\Phi^{(j1)}$ for $j = 2, \dots, n$.
Then either 

(i) $\Phi_{11}^{(21)} = 0$, or

(ii) $0 < (\Phi_{11}^{(21)})^2 \le 
\left\{\sum_{j=2}^n \sum_{k=1}^m (\Phi_{k1}^{(j1)})^2\right\}/\{(n-1)m\}
\le 1/\{(n-1)m\}$.

If (i) holds, let 
$S \in GP(m)$ such that $(\Phi^{(11)})^tS\Phi^{(11)}$ has a positive
$(1,1)$ entry equal to $\eta < \mu$. Let  $P = (P^{(ij)}) \in H$
be such that $P^{(11)} = S$, 
$P^{(22)} = \sum_{j=1}^{m/2} (e_{2j+1,2j}^{(m)} - e_{2j,  
2j+1}^{(m)})$,
$P^{(2k-1,2k)} = I_m$ and $P^{(2k,2k-1)} = -I_m$ for $k = 2, \dots, n/2$.
Then $\Phi^t P \Phi \in G$ has $(1,1)$ entry equal to $\eta < \mu$,
which contradicts the definition of $\mu$.
 
If (ii) holds and $(m,n) \ne (3,2)$, then $0 < \Phi_{11}^{(21)} <  1/2$.
Let $P = (P^{(ij)}) \in H$
be such that $P^{(12)} = I_m$, $P^{(21)} = 2e_{11}^{(m)} - I_m$,
$P^{(2k-1,2k)} = I_m$ and $P^{(2k,2k-1)} = -I_m$ for $k = 2, \dots, n/2$.
Then $\Phi^t P \Phi \in G$ has $(1,1)$ entry equal to 
$2\Phi_{11}^{(11)} \Phi_{11}^{(21)} < \Phi_{11}^{(11)} = \mu$,
which contradicts the definition of $\mu$.

If (ii) holds and $(m,n) = (3,2)$, let the first column of $\Phi$
be $v = (a_1, a_2, a_3, b_1, b_2, b_3)^t$.  By our assumption,
$\mu = a_1 \le b_1 \le b_2 \le b_3$.  Since $v$ is a unit vector,
we have $a_1 \le 1/2$.  We claim that $a_1 = 1/2$ and hence
$b_1 = b_2 = b_3 = 1/2$.  Assume that $a_1 < 1/2$. We consider
two cases.  

\noindent
Case 1.  Suppose $a_1 > 1/4$.  Note that $b_1 \ge 1/2$. Otherwise,  
one can find $R \in H$ such that the $(1,1)$ entry of $\Phi^tR\Phi$  
equals $2a_1b_1 < a_1$, which is a contradiction.  
Let $P = P_1 \oplus P_2 \in H$ with
$P_1 = -e_{11}^{(3)} - e_{23}^{(3)} + e_{32}^{(3)}$
and $P_2 = e_{11}^{(3)} - e_{23}^{(3)} + e_{32}^{(3)}$.
Then $\Phi^t P \Phi$ has $(1,1)$ entry equal to the positive number
$b_1^2 - a_1^2 \le (b_1^2 + b_2^2 + b_3^2)/3 - a_1^2 \le 
(v^tv - 4 a_1^2)/3 = (1-4a_1^2)/3 < a_1 = \mu$, 
because the roots of the equation
$4t^2 + 3t - 1 = 0$ are $-1$ and $1/4$.  

\noindent
Case 2.  Suppose $a_1 \le 1/4$.  Let $P = I_6 - 2e_{11}^{(6)}$.
Then $\Psi = \Phi^tP\Phi$ has $(1,1)$ entry equal to $1 - 2a_1^2$.
Let $u = (c_1, c_2, c_3, d_1, d_2, d_3)$ be the first column of $\Psi$.
We may assume that $u$ is a nonnegative vector.  Since
$u$ is a unit vector, there are other nonzero entries besides $c_1$.
In particular, we may assume that $c_3^2 + d_3^2 > 0$, otherwise  
multiply $\Psi$ by a suitable $P \in H$.  
Let $Q = e_{12}^{(3)} - e_{21}^{(3)} + e_{33}^{(3)}$.
Then $Q \oplus Q \in H$ and the $(1,1)$ entry of 
$\Psi^t(Q \oplus Q)\Psi$ equals $c_3^2 + d_3^2 \le u^tu - c_1^2 
= 4a_1^2(1-a_1^2) < 4a_1^2 < a_1 = \mu$.

\noindent
In both cases, we get the desired contradiction.  So, we have
$v = (1,0,0,1,1,1)^t/2$.  But then if we define $Q$ as in Case 2, the
$(1,1)$ entry of $\Phi^t(Q\oplus Q)\Phi$ will be $1/4 < 1/2 = a_1 =  
\mu$, which contradicts the definition of $\mu$.

Combining the about analysis, we get the conclusion. \qed

\medskip\noindent
{\bf Lemma 4.5} \sl Let $G$ be the isometry group of $N$ and suppose
that  $G$ is infinite or
$G < GP(mn)$.
Then every isometry in $G$ has form (\ref{new}), i.e.  Theorem~4.1$(i)$ holds.

Proof. \rm If $G$ is infinite, then by Lemma~4.3 either $G = O(mn)$ or $G =  
O(m)*P(n)$. In the former case, clearly, we have $N_1 = N_2 = \ell_2$.  In the 
latter case, we have $N_2 = \ell_2 \ne N_1$.  The conclusion of 
Theorem~4.1$(i)$ holds. 

Suppose $G < GP(mn)$.  As mentioned to the introduction to this section,
one can obtain the conclusion by the the remark after Theorem~3.2. \qed 

\medskip
Next we turn to the exceptional case when $\dim Y = m = 4$.

\medskip\noindent
{\bf Lemma 4.6} \sl Suppose ${\rm dim}\,Y = 4$. The conclusion of
Theorem 4.1$(i)$ holds $($i.e. isometries have form $(9)\,)$.

Proof. \rm Let $G$ be the isometry group of $N$. 
By Lemma~4.5 it remains to consider the case
when $G$ is finite and is not a subgroup of $GP(mn)$.  

For that, define $\mu$ 
as in the proof of Lemma 4.4, and let $\Phi \in G$ have $(1,1)$ entry
equal to $\mu$ and a nonnegative first column. We divide the proof into 
three assertions.  The matrix
$$Q = e_{12}^{(4)} + e_{34}^{(4)} -  e_{21}^{(4)} - e_{43}^{(4)}$$  
will be used frequently in our arguments.

\medskip\noindent
{\sl Assertion} 1. 
The first column of $\Phi$ equals $v = (1,1,1,1,0,0,\dots,)^t/2$.
In particular, $\mu = 1/2$.

Suppose the first column of $\Phi$ equals $v = (v_1, v_2, \dots)^t$
with $0 \le v_2 \le v_3 \le v_4$.  Then $v_2 \ne 0$.  Otherwise,
we may let $P = (e_{34}^{(4)} - e_{43}^{(4)}) \oplus (I_{n-1}\otimes Q)
\in H$, where $Q$ is defined as above, so 
that the $(1,1)$ entry of $\Phi^t P \Phi$ equals $v_1^2 < v_1 = \mu$,
which contradicts minimality of $\mu$.

Thus, we may assume that $v_1 \le v_2 \le v_3 \le v_4$.  

To prove that $v$ is of the asserted form it is enough to show that
$v_1 = 1/2$, since $v^tv = 1$.

Since $v_1^2 \le (v_1^2 + v_2^2 + v_3^2 + v_4^2)/4 = v^tv/4 = 1/4$, we see 
that $v_1 \le 1/2$. One can show that $v_1 < 1/2$ is impossible as in the last 
part (Case~2) of the proof of Lemma 4.4 (see also the first proof of Theorem 
3.2 in \cite{DLR}). 

\medskip\noindent
{\sl Assertion} 2. A column of $\Psi \in G$ must be of the form
$e_i^X\otimes (Pu)$ for some $1 \le i \le n$, $P \in GP(4)$ and 
$u = (1,0,0,0)^t, (1,1,1,1)^t/2$, or $(1,1,0,0)^t/\sqrt 2$.

\medskip
Let $v$ be the $k$th column of $\Psi \in G$. Multiplying $\Psi$ by a suitable
$R \in H$, we may assume that $k = 1$ and the $(1,1)$ entry of $\Psi$ is 
nonzero.  We need to show that $v = e_1^X \otimes (Px)$.
If $v$ is not of this form, then (see the first proof of Theorem 3.2 in 
\cite{DLR}) there exists $S \in GP(4)$ such that the $(1,1)$ entry of 
$(\Psi^{(11)})^tS \Psi^{(11)}$ 
equals $\eta < 1/2$.  Let $R = S \oplus (I_{n-1} \otimes Q) \in H$,  where
$Q$ is defined as above.
Then the $(1,1)$ entry of $\Psi^t R \Psi \in G$ equals $\eta < 1/2$,
which is a contradiction.  

\medskip\noindent
{\sl Assertion} 3. Suppose $\Psi = (\Psi^{(ij)}) \in G$ is in
$n\times n$ block form with $\Psi^{(ij)} \in \IR^{4\times 4}$
for each $(i,j)$.  Then each nonzero $\Psi^{(ij)}$ must be of the  
form
$P,  PAR$ or $PBR$, with $P, R \in GP(4)$, where
$$A = (I_4 - (1,1,1,1)^t (1,1,1,1)/2)  
\quad\hbox{ and } \quad 
B = {1\over 2}\left\{\pmatrix {1 & 1 \cr 1 & -1 \cr} 
\oplus \pmatrix {1 & 1 \cr 1 & -1 \cr}
\right\}.$$

Consider a nonzero $\Psi^{(ij)}$.  Multiplying $\Psi$ by a suitable
$R \in H$, we may assume that $(i,j) = (1,1)$, the first column of 
$\Psi$ is nonnegative, and the $(1,1)$ entry $\alpha$
of $\Psi$ has the smallest
magnitude among all nonzero entries in $\Psi^{(11)}$. We consider 3  
cases.

\noindent
Case 1.  Suppose $\alpha = 1/2$.  Then the first column of 
$\Psi$ equals $(1,1,1,1,0,\dots)^t/2$. Since $\Psi^t \in G$ as well,
we see that the first columns of $\Psi^t$ are of the form
$e_1^X \otimes P(1,1,1,1)^t/2$ for some $P \in GP(4)$ by
Assertion 2.  Thus $\Psi^{(11)} = PAR$ for some $P, R \in GP(4)$.

Note that it follows from Assertion 1 and the above arguments
that $\Phi^{(11)}$ is of the form $PAR$ for some $P, R \in GP(4)$.

\noindent
Case 2. Suppose $\alpha = 1/\sqrt 2$. By Assertion 2, we may assume  
that the first column of $\Psi^{(11)} = (1,1,0,0)^t/\sqrt 2$.  Since 
$\Psi^t \in G$ as well, we may assume the first row of $\Psi^{(11)}$ 
equals $(1,1,0,0)/\sqrt 2$.  Since the first two column of $\Psi$ are
orthogonal, the second column of $\Psi^{(11)}$ equals
$(1,-1,0,0)^t/\sqrt 2$.  By Assertion 2 and the knowledge about 
the first two columns of $\Psi$, one easily sees that the $(1,1)$ and $(2,1)$  
entries of $\Psi^t \Phi \in G$ are $1/\sqrt 2$ and $0$.  By Assertion 2  
again,  one of the $(3,1)$ and $(4,1)$ entries is $0$, and the other
has magnitude $1/\sqrt 2$. It follows that the third and the four columns
of $\Psi$ have the same form. Thus $\Psi^{(11)} = PBR$ for some
$P, R \in GP(4)$.  

\noindent
Case 3. Suppose $\alpha = 1$. Then the first column of $\Phi\Psi \in G$
equals $(1,1,1,1,0,\dots)^t/2$.  By the result in Case 1, the first  
four columns of $\Phi\Psi$ are of the form $e_1^X \otimes P(1,1,1,1)/2$
for some $P \in GP(4)$.  Thus $\Psi_{(11)} \in GP(4)$.

\medskip
We are now ready to complete the proof of the lemma.
Let ${\cal A}$ be the group generated by $GP(4)$ and $A$, and let
${\cal B}$ be the group generated by $GP(4)$ and $B$.  It is known
(e.g., see \cite[Theorem 3.2]{DLR}) that ${\cal A}$ is a normal subgroup
of ${\cal B}$, and they are the only other possible isometry groups of
a symmetric norm on $\IR^4$ besides $O(4)$ and $GP(4)$.  
By Assertion 3, one easily concludes
that if $G$ is not infinite and $G$ is not a subgroup of $GP(mn)$,
then $G$ must be of the form $G_2*P(n)$,
where $G_2 = GP(4), {\cal A}$ or ${\cal B}$.  In each case, 
$G_2$ is clearly the isometry group of $N_2$. \qed

\medskip
Finally, we deal with the exceptional case when $\dim Y = m = 2$ in 
the next two lemmas.

\medskip\noindent
{\bf Lemma 4.7} \sl Suppose ${\rm dim}\,Y = 2$. Then the isometry group of $N$ has form (\ref{new}) or (\ref{new2}).

Proof. \rm  Let $G$ be the isometry group of $N$.  
By Lemma 4.5, we only need to consider the case when $G$ 
is finite and not a subgroup of $GP(mn)$.
Define $\mu$ as in the proof of Lemma~4.4. Then $0 < \mu < 1$.  Let
$\Phi \in G$ have $(1,1)$ entry equal to $\mu$ and a nonnegative first  
column $v = (v_1, \dots, v_{2n})^t$.  Similarly as in Lemma~4.6,
we divide the proof into several assertions.

\smallskip\noindent
{\sl Assertion} 1.  The vector $v$ cannot have more than four nonzero  
entries. 

If the assertion is not true, then
$$\eta = \min\{v_j: 2 \le j \le 2n, v_j > 0\} \le \left\{(\sum_{j=2}^{2n} 
v_j^2)/(k-1)\right\}^{1/2} < 1/2,$$ 
where $k$ is the number of nonzero entries of $v$.  
Let $H = GP(2)*GP(n)$.  If $v_2 = \eta$, we can
find $P \in  H$ such that the first column of $P\Phi$
equals $(v_2, v_1, v_4, -v_3, v_6, -v_5, \dots, v_{2n},  -v_{2n-1})^t$.
If $v_2 > \eta = v_j$ for some $j > 2$, we may assume that $j = 3$  
after multiplying $\Phi$ by a suitable $Q \in H$ on the left. Then we can  
find  $P \in H$ such that the first column of $P\Phi$ equals
$(v_3, -v_4, v_1, v_2, v_6, -v_5, \dots, v_{2n})^t$.  In both cases, 
$\Phi^t P \Phi$ has $(1,1)$ entry equal to $2\mu\eta < \mu$, 
which  contradicts minimality of $\mu$.                                    

\smallskip\noindent
{\sl Assertion} 2. We have $v_2 \ge 1/2$.  

If $v_2 = 0$, we can 
find $P \in  H$ such that the first column of $P\Phi$ equals 
$$(v_1, v_2, v_4, -v_3, v_6, -v_5, \dots, v_{2n}, -v_{2n-1})^t$$
so that the $(1,1)$ entry of $\Phi^tP\Phi$ is $\mu^2 < \mu$, 
which is a contradiction.  If $0 < v_2 < 1/2$, we can 
find $P \in  H$ such that the first column of $P\Phi$
equals $(v_2, v_1, v_4, -v_3, v_6, -v_5, \dots, v_{2n}, -v_{2n-1})^t$  
so that $\Phi^tP\Phi$ has $(1,1)$ entry equal to $2v_2\mu < \mu$, 
which is a contradiction.

\smallskip\noindent
{\sl Assertion} 3. 
The vector $v$ cannot have exactly 3 nonzero entries.

If the assertion is not true,  we may assume without loss of generality
 that $v_j = 0$ for 
$j = 4, \dots, 2n$(by replacing $\Phi$ with $P\Phi$  for some 
$P  \in H$, if necessary).

If $n \ge 3$, let $R \in H$ be such that the first column of $R\Phi$ equals
$(0, 0, v_1, v_2, v_3, 0, \dots)^t$, then the $(1,1)$ entry of $\Phi^tR\Phi$ 
equals $\mu v_3 < \mu$, which is a contradiction. 

If $n = 2$, let $S \in H$ be such that the first column of 
$S\Phi$ equals $(-v_1, v_2, v_3, 0)^t$.  Then one can find $R \in H$
such that the first column of $\Psi = R\Phi^tS\Phi \in G$ is nonnegative 
and equals $(1-2\mu^2, r, s, t)$ with $s \ge t$.  Furthermore, 
suppose $\eta = v_2 r + v_3 s - \mu(1-2\mu^2)$ is the $(1,1)$
entry of $\Psi^tS \Phi \in G$.  Since 
$$(v_2r+v_3s)^2 \le (v_2^2+v_3^2)(r^2+s^2)\quad \hbox{ and }\quad
r^2+s^2+t^2+(1-2\mu^2)^2 = 1 = v_1^2 + v_2^2 + v_3^2,$$
we have
\begin{eqnarray} 
\eta
& \le & \{(v_2^2 + v_3^2)(r^2 + s^2)\}^{1/2} - \mu(1-2\mu^2)\cr
& \le & \{(1-v_1^2)(1-(1-2\mu^2)^2)\}^{1/2} - \mu(1-2\mu^2) = \mu,  
\label{ck}
\end{eqnarray}
and $\eta = \mu$ if and only if $(r,s,t)$ is a multiple of
$(v_2, v_3, 0)$.

Note that $v_3 \ge 1/2$. Indeed, otherwise, one can 
find $P \in  H$ such that the first column of $P\Phi$
equals $(v_3, 0, v_1, v_2)^t$ so 
that $\Phi^tP\Phi$ has $(1,1)$ entry equal to $2v_3\mu < \mu$, 
which is a contradiction.

If $s = 0$, then $t = 0$ and  $r^2 = 1 - (1-2\mu^2)^2 = 4\mu^2(1-\mu^2)$.
Then $(r,s,t)$ is not a multiple of $(v_2,v_3,0)$ and we have
$\mu > \eta = v_2 r -\mu(1-2\mu^2)\ge r/2 - \mu(1-2\mu^2) 
= \mu\sqrt{1-\mu^2} -\mu(1-2\mu^2) > 0$, since $\mu\le 1/\sqrt{3}$,
which  contradicts minimality of $\mu$.

If $r = 0$, then $(r,s,t)$ is not a multiple of $(v_2, v_3, 0)$,
and hence $\mu > \eta$.  If $\eta = v_3 s - \mu(1-2\mu^2) \ne 0$, 
then we can find $P \in H$ such that the $(1,1)$ entry of 
$\Psi^tP\Phi$ equals $|\eta| < \mu$, which is a contradiction.
Similarly, if $v_3 t - \mu(1-2\mu^2) \ne 0$, we have a contradiction.
Suppose $t^2 = s^2 = 2\mu^2(1-\mu^2)$.  Then 
$v_3^2 =  (1-2\mu^2)^2/(2(1-\mu^2))$.  Since $v_3^2  - \mu^2 \ge 0$,
we have $6\mu^4 - 6 \mu^2 + 1 \ge 0$. Thus $\mu^2 \le 1/2 - 1/\sqrt{12}$
or $\mu^2 \ge 1/2 + \sqrt{12}$.  Since $\mu^2 \le (v_1^2 + v_2^2 + v_3^2)/3 
= 1/3$, we have $\mu^2 \le 1/2 - 1/\sqrt{12} \le 1/4$.  Now, if $\mu < 1/4$
or $1/4 \le \mu < 1/2$, one can derive a contradiction as in the last
part of the proof of Lemma 4.4 (see also the first proof of Theorem 3.2 in 
\cite{DLR}).

If $r, s \ne 0$, then $r, s \ge \mu$. Since $v_2, v_3 \ge 1/2$, we have
$$\eta \ge \mu(v_2+v_3) - \mu(1-2\mu^2) \ge \mu - \mu(1-2\mu^2) >  
0.$$
If $(r,s,t)$ is not a multiple of $(v_2,v_3,0)$, then, by (\ref{ck}), 
$\mu > \eta > 0$, which is a contradiction.  Thus $(r, s, t)$ is 
a multiple of $(v_2, v_3, 0)$ and $t = 0$.
We can find   $P \in H$ such that the $(1,1)$ entry of 
$\Psi^tP\Phi$ equals $\delta = |v_2r - \mu(1-2\mu^2)| < \mu$, which 
is a contradiction if $\delta > 0$.  If $\delta = 0$, one can let
$(r,s) = c(v_2, v_3)$ and solve the equations: 
$$v_2^2 + v_3^2 = 1 - \mu^2$$
$$r^2 + s^2 = c^2(v_2^2+v_3^2) = 1 - (1-2\mu^2)^2 = 4  
\mu^2(1-\mu^2),$$
$$0 = v_2r - \mu(1-2\mu^2) = cv_2^2 - \mu(1-2\mu^2)$$ 
to conclude that $v_2^2 = 1/2 - \mu^2$ and $v_3^2 = 1/2$.   Since, by 
Assertion~2,  $v_2 \ge \mu$, we have $\mu \le 1/2$ and since 
$2\mu^2 +1/2\le v_1^2 + v_2^2 + v_3^2 = 1$ 
we conclude that $\mu = 1/2 = v_2$. 

Now, we may assume that the first row of $\Phi$ is $v^t$.  Otherwise
replace $\Phi$ by $\Phi R$ for some $R \in H$.  One readily sees that
$P\Phi Q$ equals $L_1$ or $L_2$ for some $P, Q \in H$, where   
$$L_1 = 2^{-1}\pmatrix{1 & 1 & \sqrt 2 & 0 \cr
                      1 & -1 & 0 & \sqrt 2 \cr
                     \sqrt 2 & 0 & -1 & -1 \cr
                     0 & \sqrt 2 & -1 & 1 \cr}  \quad \hbox{ and }  
\quad
L_2 = 2^{-1}\pmatrix{1 & 1 & \sqrt 2 & 0 \cr
                    1 & 1 & -\sqrt 2 & 0 \cr
              \sqrt 2 & -\sqrt 2 & 0 & 0 \cr
                           0 & 0 & 0 & 1 \cr}.
$$
It is routine to check that $\langle H, L_1\rangle = \langle H, L_2 \rangle$
consists of matrices of the form $P, PL_1Q$ or $PL_2 Q$ with $P, Q \in H$.  
By the previous discussion, this group is contained in $G$.
We shall show that this is impossible.  Suppose 
$\langle H, L_1\rangle = \langle H, L_2 \rangle < G$.  Then
\begin{eqnarray*}
N_2(e_1^Y + e_2^Y) 
& = & N(e_{11} - e_{21}) = N(L_2(e_{11} - e_{21})) \cr
& = & N(\sqrt 2 e_{12}) = \sqrt 2 N_2(e_1^X) = \sqrt 2 .\cr
\end{eqnarray*}
Thus 
\begin{eqnarray*}
N_1(ae_1^X+be_2^X) 
& = & N(ae_{11} + be_{12}) \cr
& = & N(2^{-1/2}a(e_{11} - e_{21}) + b e_{22}) \cr
& = & N(L_2(2^{-1/2}a(e_{11} - e_{21}) + b e_{22})) \cr
& = & N(ae_{12} + be_{22}) = N_2(ae_1^Y+be_2^Y) \cr
\end{eqnarray*}
for any $a,b \in \IR$. Furthermore, 
\begin{eqnarray*}
N_2(ae_1^Y+be_2^Y) 
& = & N(ae_{11}+be_{21}) \cr
& = & N(L_2(e_{11}+be_{21})) \cr
& = & N((a+b)(e_{11}+e_{21})/2 + (a-b)e_{12}/\sqrt 2) \cr
& = & N((a+b)e_{11} + (a-b)e_{12})/\sqrt 2 \cr
& = & N_1((a+b)e_1^X + (a-b)e_2^X)/\sqrt 2 \cr
& = & N_2((a+b)e_1^X + (a-b)e_2^X)/\sqrt 2 \cr
\end{eqnarray*} 
for any $a, b \in \IR$.
Thus $A = 2^{-1/2} \pmatrix{ 1 & 1 \cr 1 & -1 \cr}$ is an isometry for
$N_2$, and $\Gamma = A \oplus I_2$ is an isometry for $N$. Note that
the first column of $\Gamma$ is neither of the form $P(1,0,\dots)^t$
nor $P(1/2,1/2,1/\sqrt 2, 0, \dots)^t$ with $P \in H$.
We can find $R \in H$ such that the $(1,1)$ entry of 
$\Gamma L_2$ is positive and less than $1/2$, which is a  
contradiction with the fact that $\mu = 1/2$.

\medskip\noindent
{\sl Assertion} 4. Suppose the first column of $\Phi$ has exactly 4  
nonzero entries. Then the first column of $P\Phi$ equals
$(1,1,1,1,0,\dots)^t/2$ for some $P \in H$, and every $\Psi \in G$ has form (\ref{new2}).

Suppose $v$ has exactly 4 nonzero entries.  One can show that all of them 
equal $1/2$ by arguments similar to those in the proof of Lemma 4.6  
(cf.\ Assertion 1). By Assertion 2, we have $v_1 = v_2 = 1/2$.  Since
$\mu = 1/2$, by Assertion 2,  for each  $1 \le j \le n$, 
$v_{2j-1}$ and $v_{2j}$ are either both zero or both nonzero. Thus 
$P\Phi$ has first column equal to $(1,1,1,1, 0,\dots)^t/2$ for some  
$P \in H$ as asserted.

Since $(P\Phi)^t \in G$ and the first four entries in
the first row of $(P\Phi)^t$ equal $1/2 = \mu$, the second column $w$
of $P\Phi$ is of the form
$(\pm 1, \pm 1, \pm 1, \pm 1, 0, \dots)^t/2$ by Assertion 2.
We claim that $w = \pm(1,1,-1,-1, 0, \dots)^t/2$.  If it is not
true, we may assume that $w = (1, -1, 1, -1, 0, \dots)^t/2$.
Otherwise interchange the third and fourth rows of $P \Phi$, and
multiply the second column of $P \Phi$ by $-1$ if necessary.
It follows that 
\begin{eqnarray*} N_2((a,b)^t) 
& = & N(ae_{11} + be_{21}) = N(P\Phi(ae_{11} + be_{21})) \cr
& = & N(a+b)(e_{11} + e_{12}) + (a-b)(e_{21} + e_{22}))/2 \cr
& = & N_2((a+b, a-b)^t) N_1(1,1)/2. \cr
\end{eqnarray*}
It follows that $A = k\pmatrix{1&1\cr 1&-1\cr}$ is an isometry for 
$N_2$ for $k=N_1((1,1)^t)/\sqrt 2$.  Since an isometry for $N_2$ must  
be orthogonal (e.g., see \cite{DLR}), $k = 1/\sqrt 2$.  But then $\Psi =  
A \oplus I_{2n-2} \in G$ and the first column of $\Gamma = \Psi P \Phi$ has 
3  nonzero  entries, and there exists $Q \in H$ such that $\Gamma^t Q \Phi$ 
has a positive nonzero entry less than $1/2 = \mu$, which is impossible. 

Now for any $\Psi \in G$, the columns of $\Psi$ must be of the form
$Q(1,1,1,1,0,\dots)^t/\sqrt 2$, $Q(1,1,0, \dots)^t/\sqrt 2$,
or $Q(1,0,\dots)^t$ for some $Q \in H$.  Otherwise, one can
find $R, S \in H$ such that $R\Psi^t S\Phi$ has positive $(1,1)$
entry less than $1/2 = \mu$.  Now if $\Psi$ has a column of the form
$Q(1,1,0,\dots)^t/\sqrt 2$, then one can show 
that $2^{-1/2}\pmatrix{1&1\cr 1&-1\cr}$ is an isometry for
$N_2$, and get a contradiction as in the preceding paragraph.
Furthermore, if the $(2j-1)$th (respectively, $(2j)$th) column of  
$\Psi$ is of the form $Q(1,1,1,1,0, \dots)^t/2$, then the $(2j)$th  
(respectively, $(2j-1)$th column must be
of the form $\pm Q(1,1,-1,-1, 0, \dots)^t/2$ by arguments similar to 
those in the analysis of the second column of $\Phi$.
It follows that every $\Psi \in G$ permutes the matrices
in the set $\{\pm(e_1^Y\pm e_2^Y)(e_j^X)^t: 1 \le j \le n\}$ as asserted.

\smallskip\noindent
{\sl Assertion} 5.  Suppose $v$ has exactly 2 nonzero entries.
The conclusion of Theorem 4.1$(i)$ holds.

By Assertion 2, if $v$ has exactly two nonzero entries, then 
we may assume that $(v_1, v_2) = (\sin t, \cos t)$ for some 
$t \in (0, \pi/4)$.  Now, it is easy to see the columns of $\Psi \in  G$
must be of the form $P(a,b, 0, \dots)^t$ with $a^2 + b^2 = 1$.
Otherwise, one can find $R, S \in H$ such that the $(1,1)$
entry of $R\Psi^tS\Phi$ is positive and is less than $\sin t$.
Moreover, if the $(2j-1)$th (respectively, the $(2j)$th)
column of $\Psi$ is of the form $P(a,b, 0, \dots)^t$ with $ab \ne 0$,
then by the fact that $\Psi \in G$ one can conclude that the
$(2j)$th (respectively, the $(2j-1)$th) column must be of the form
$\pm P(b,-a, 0, \dots)^t$. Thus $P\Psi P^t$ is a direct sum of 
a signed permutation matrix $A$ and a number of $2\times 2$ orthogonal 
matrices $B_i$.  Furthermore, we may assume that $A$ is a direct sum
of matrices in $GP(2)$.  Otherwise, $\Gamma = A \oplus I_k \in G$
will satisfy the hypothesis of Lemma 2.4, and hence both 
$N_1$ and $N_2$ equal $\ell_p$ for some $p \ge 1$.  Thus $P\Psi P^t$ 
must be a direct sum of isometries for $Y$ for some $P \in H$,
and the conclusion follows.
\qed

\medskip\noindent
{\bf Lemma 4.8} \sl If isometries of $X(Y)$ have form (\ref{new2}), then
$X = \ell_p^n$ and $Y = E_p(2)$ for some $p, \ \ 1\le p<\infty, \ \ p\neq2.$

Proof. \rm  Suppose $m = 2$.  
The elements in $Y$ will be written as $(y_1,y_2)$, and
the elements in $X(Y)$ will be written as
$$(x_{11}, x_{21}, x_{12}, x_{22}, \dots, x_{1n}, x_{2n}).$$
If isometries of $X(Y)$ have form (\ref{new2}) but condition (i) of Theorem~4.2 does not hold, then the isometry group 
$G$ of $X(Y)$ must be $A_{2n}$ as mentioned before Lemma~4.2.
Therefore, the linear map 
$T: X(Y) \lra X(Y)$ defined by
$$\matrix{
& ~~~  T(x_{11},x_{21},x_{12},x_{22},\dots,x_{1n},x_{2n}) \hfill \cr  & \cr
& =  
\biggl(2^{-1}(x_{11}+x_{21}+x_{12}+x_{22}),
 2^{-1}(x_{11}+x_{21}-x_{12}-x_{22}),  \hfill\cr
& ~~~~ 2^{-1}(x_{11}-x_{21}+x_{12}-x_{22}), 
 2^{-1}(x_{11}-x_{21}-x_{12}+x_{22}), 
 x_{13},x_{23},\dots,x_{1n},x_{2n}) \biggr) \cr
}$$
is an isometry of $X(Y)$. 
Assume that $x_{13}=x_{23}=\dots=x_{1n}=x_{2n}=0$, then
$$N(x_{11},x_{21},x_{12},x_{22},0,0,\dots,0,0) 
= N_1(N_2(x_{11},x_{21})e_1^X + N_2(x_{12},x_{22})e_2^X).$$
In particular, when $x_{11}=x_{21}= a$ and $x_{12}=x_{22}=b$ we get
$$N(a,a,b,b,0,\dots,0) 
= N_1(aN_2(1,1)e_1^X + bN_2(1,1)e_2^X)
= N_2(1,1)N_1(a e_1^X + b e_2^X).$$
Since $T$ is an isometry
\begin{eqnarray*}
N(a,a,b,b,0,\dots,0)  &=& N(T(a,a,b,b,0,\dots,0)) \\
&=& N(a+b,a-b,0,\dots,0)= N_1(N_2(a+b,a-b)e_1^X) \\
&=& N_2(a+b,a-b).
\end{eqnarray*}
Thus
\begin{equation} \label{EX1}
N_1(a e_1^X + b e_2^X) = \frac1{N_2(1,1)}N_2(a+b,a-b)
\end{equation}
for all $a,b \in\bbR.$ Also
\begin{equation} \label{EX}
N_2(c,d) = N_2(1,1)N_1(\frac12(c+d) e_1^X + \frac12(c-d) e_2^X)
\end{equation}
for all $c,d \in\bbR.$
Further, since $T$ is an isometry, we get
$$ N(x_{11},x_{21},x_{12},x_{22},0,0,\dots,0,0) 
= N(T(x_{11},x_{21},x_{12},x_{22},0,0,\dots,0,0)).$$
It follows that
\begin{eqnarray} \label{iso}
& & ~~~~~N_1\biggl(N_2(x_{11}, x_{21})e_1^X + 
N_2(x_{12},x_{22})e_2^X\biggr) \nonumber \\
& &  =
N_1\biggl[N_2\biggl(2^{-1}(x_{11}+x_{21}+x_{12}+x_{22}),
2^{-1}(x_{11}+x_{21}-x_{12}-x_{22})\biggr)e_1^X \nonumber \\
& &  \quad \quad \quad + N_2\biggl(2^{-1}(x_{11}-x_{21}+x_{12}-x_{22}),
   2^{-1}(x_{11}-x_{12}-x_{21}+x_{22}\biggr)e_2^X\biggr]
\end{eqnarray}
Put
\begin{eqnarray*}
x_{11}+x_{21} &=& y_1 \quad \quad \qquad x_{11}-x_{21} =  y_2 \cr
x_{12}+x_{22} &=& y_3 \quad \quad \qquad \ x_{12}-x_{22} = y_4.
\end{eqnarray*}
Then by (\ref{EX}) and (\ref{iso}) we get
\begin{eqnarray} \label{iso2}
& & ~~~~2^{-1}N_2(1,1)N_1\biggl(N_1(y_1e_1^X+y_2e_2^X)e_1^X + 
      N_1(y_3e_1^X+y_4e_2^X)e_2^X \biggr) \nonumber \\
& &  = 
2^{-1}N_2(1,1)N_1\biggl(N_1(y_1e_1^X+y_3e_2^X)e_1^X + 
      N_1(y_2e_1^X+y_4e_2^X)e_2^X \biggr)
\end{eqnarray}
Define $f(a,b) = N_1(a e_1^X + b e_2^X).$ Then (\ref{iso2}) becomes:
$$ f(f(y_1,y_2),f(y_3,y_4)) = f(f(y_1,y_3),f(y_2,y_4))$$
for all $y_1,y_2,y_3,y_4\in \bbR.$ When $y_4 = 0 $   we have
$$f(f(y_1,y_2), y_3 )= f(f(y_1,y_2),f(0,y_3)) = f(f(y_1,0),f(y_2,y_3)) =
 f(y_1,f(y_2,y_3)).$$
By a theorem of Bohnenblust [Bo]
$$f(a,b) = \max(|a|,|b|) {\rm \ \ \ \ 
or\ \ \ \ } f(a,b) =  (|a|^p+|b|^p)^{1/p}$$
for some $p, \ \ 1\le p<\infty.$
Hence 
\begin{equation} \label{LP}
N_1(a e_1^X + b e_2^X)  =  \ell_p(a,b)
\end{equation}
for some $p, \ \ 1\le p\le\infty.$
By (\ref{EX}) we see that $Y=E_p(2).$ 

To see that $X=\ell_p^n$, let $k\le n$ be the maximal number such that
\begin{equation} \label{AS}
N_1(a_1 e_1^X +\dots+ a_k e_k^X)  =  \ell_p(a_1,\dots,a_k)
\end{equation}
for all $a_1,\dots,a_k\in \bbR.$
If $k=n$ there is nothing to prove. If $k<n$ let 
$(a_l)_{l=1}^{k+1}\subset\bbR^{k+1}$ be arbitrary and set 
$x_{1l}=x_{2l}=a_l$ for $l=1,\dots,k+1,$ and $x_{1l}=x_{2l}=0$ 
for $l>k+1.$ Then
\begin{eqnarray} \label{T}
&~~ & N(a_1,a_1,\dots,a_{k+1},a_{k+1},0,\dots,0) \nonumber \\
&=& N_1(a_1N_2(1,1) e_1^X +
   \cdots+ a_{k+1}N_2(1,1) e_{k+1}^X) \nonumber \\
&=& N_2(1,1)N_1(a_1 e_1^X +\cdots+ a_{k+1} e_{k+1}^X) 
\end{eqnarray}
Since $T$ is an isometry we get
\begin{eqnarray*}  
&~~ & N(a_1,a_1,\dots,a_{k+1},a_{k+1},0,\dots,0)\\
&=& N(T(a_1,a_1,\dots,a_{k+1},a_{k+1},0,\dots,0)) \\
&=& N((a_1+a_2,a_1-a_2,0,0,a_3,a_3,\cdots,a_{k+1},a_{k+1},0,\dots,0)) \\
&=& N_1(N_2(a_1+a_2,a_1-a_2) e_1^X +
    a_3N_2(1,1) e_3^X+\cdots+ a_{k+1}N_2(1,1) e_{k+1}^X) \\
&=& N_1(N_2(1,1)\ell_p(a_1,a_2) e_1^X +a_3N_2(1,1) e_3^X+ \cdots + 
    a_{k+1}N_2(1,1) e_{k+1}^X) \qquad   
     {\rm by \  (\ref{EX1})\   and\   (\ref{LP})}\\
&=& N_2(1,1) N_1(\ell_p(a_1,a_2) e_1^X +a_3  e_3^X+\cdots+ 
    a_{k+1}  e_{k+1}^X) \\
&=& N_2(1,1) \ell_p(\ell_p(a_1,a_2),a_3,\dots,a_{k+1}) \hskip 1.7in
    {\rm by \ (\ref{AS})\   and \   symmetry \  of\    }  X  \\
&=& N_2(1,1) \ell_p(a_1,a_2,a_3,\dots,a_{k+1}). 
\end{eqnarray*}
Hence, by (\ref{T})
$$N_1(a_1 e_1^X +\dots+ a_{k+1} e_{k+1}^X) = \ell_p(a_1,\dots,a_{k+1}),$$
which contradicts maximality of $k$.
Thus $(X,N_1)=\ell_p^n.$

Finally, we conclude that $p\neq 2.$ Indeed, if $p=2$ then 
$E_p(2) = \ell_2^2$ and then $X(Y)= \ell_2(\ell_2^2) = \ell_2^{2n}$
whose isometry group  is not of the form (\ref{new2}).\qed

\bigskip \noindent
{Department of Mathematics, The College of William and Mary,
Williamsburg, VA 23187 \\
E-mail: ckli@cs.wm.edu}

\medskip\noindent
{Department of Mathematics, The University of Texas at Austin, 
Austin, TX 78712 \\
E-mail: brandri@math.utexas.edu}

\end{document}